\documentclass[11pt,a4paper]{amsart}
\usepackage{amsthm,amsmath,amsfonts,float,calc}
\usepackage[lmargin=35mm,rmargin=35mm,tmargin=29mm,bmargin=29mm]{geometry}
\usepackage[sort&compress,numbers]{natbib}
\renewcommand{\baselinestretch}{1.2}
\theoremstyle{plain}
\newtheorem{theorem}{Theorem}
\newtheorem{lemma}{Lemma}
\newtheorem{corollary}{Corollary}
\newtheorem{example}{Example}
\newtheorem{proposition}{Proposition}

\theoremstyle{definition}
\newtheorem{definition}{Definition}

\newcommand{\thmlabel}[1]{\label{thm:#1}}
\newcommand{\thmref}[1]{Theorem~\ref{thm:#1}}
\newcommand{\twothmref}[2]{Theorems~\ref{thm:#1} and \ref{thm:#2}}
\newcommand{\lemlabel}[1]{\label{lem:#1}}
\newcommand{\lemref}[1]{Lemma~\ref{lem:#1}}
\newcommand{\twolemref}[2]{Lemmas~\ref{lem:#1} and \ref{lem:#2}}
\newcommand{\eqnlabel}[1]{\label{eqn:#1}}
\newcommand{\Eqnref}[1]{Equation~\eqref{eqn:#1}}
\newcommand{\twoEqnref}[2]{Equations~\eqref{eqn:#1} and \eqref{eqn:#2}}
\newcommand{\seclabel}[1]{\label{sec:#1}}
\newcommand{\secref}[1]{Section~\ref{sec:#1}}
\newcommand{\mySection}[2]{\section{#1}\seclabel{#2}}

\newcommand{\corlabel}[1]{\label{cor:#1}}
\newcommand{\corref}[1]{Corollary~\ref{cor:#1}}
\newcommand{\tablabel}[1]{\label{tab:#1}}
\newcommand{\tabref}[1]{Table~\ref{tab:#1}}
\newcommand{\proplabel}[1]{\label{prop:#1}}
\newcommand{\propref}[1]{Proposition~\ref{prop:#1}}
\newcommand{\exlabel}[1]{\label{ex:#1}}
\newcommand{\exref}[1]{Example~\ref{ex:#1}}

\newcommand{\PrimeProduct}[1]{\ensuremath{\protect\Pi_{#1}}}
\newcommand{\CARD}[1]{\ensuremath{\protect\left|#1\right|}}
\newcommand{\ceil}[1]{\ensuremath{\protect\lceil#1\rceil}}

\newcommand{\floor}[1]{\ensuremath{\protect\lfloor#1\rfloor}}
\newcommand{\e}{\ensuremath{\boldsymbol{e}}}
\newcommand{\paran}[1]{\textup{(}#1\textup{)}}
\newcommand{\bracket}[1]{\ensuremath{\protect\left(#1\right)}}
\newcommand{\BRACKET}[1]{\ensuremath{\protect\left(#1\right)}}
\newcommand{\Oh}[1]{\ensuremath{\protect\mathcal{O}(#1)}}
\newcommand{\CartProd}{\ensuremath{\,\square\,}}
\newcommand{\CP}[2]{\ensuremath{#1{\CartProd}#2}}
\newcommand{\CCP}[3]{\ensuremath{#1{\CartProd}#2{\CartProd}\cdots{\CartProd}#3}}
\newcommand{\uu}{\ensuremath{\tilde{u}}}
\newcommand{\vv}{\ensuremath{\tilde{v}}}
\newcommand{\ww}{\ensuremath{\tilde{w}}}
\newcommand{\xx}{\ensuremath{\tilde{x}}}

\newcommand{\GG}{\ensuremath{\tilde{G}}}
\newcommand{\F}{\ensuremath{\mathcal{F}}}
\newcommand{\C}{\ensuremath{\mathcal{C}}}
\newcommand{\PP}{\ensuremath{\mathbb{P}}}
\newcommand{\Density}{\delta}
\newcommand{\UpperDensity}{\overline{\delta}}
\newcommand{\N}{\ensuremath{\mathbb{N}}}
\newcommand{\Z}{\ensuremath{\mathbb{Z}}}
\newcommand{\NZ}{\ensuremath{\mathbb{N}_0}}

\newcommand{\EndProof}[1]{
\begin{minipage}[b]{\textwidth-1cm}
#1
\end{minipage}\hfill\qed
\vspace*{1ex}
\renewcommand{\qed}{}
}

\begin{document}

\title{Colourings of the Cartesian Product of Graphs\\ and Multiplicative Sidon Sets}

\author{Attila P\'or}
\address{\newline Department of Mathematics\newline  Western Kentucky University\newline Bowling Green,  Kentucky, U.S.A.}
\email{attila.por@wku.edu}

\author{David R.~Wood}
\address{\newline QEII Research Fellow\newline Department of Mathematics and Statistics\newline The University of Melbourne\newline Melbourne, Victoria, Australia}
\email{D.Wood@ms.unimelb.edu.au}

\date{November 10, 2005. Revised: \today}

\keywords{graph, colouring, cartesian product, distance-$2$ colouring, acyclic colouring, star colouring, L$(p,1)$-labelling, grid, toroidal grid, multiplicative Sidon set}


\thanks{This is an amplified version of a paper to appear in \emph{Combinatorica}. Research initiated at the Department of Applied Mathematics and the Institute for Theoretical Computer Science, Charles University, Prague, Czech Republic. Supported by project LN00A056 of the Ministry of Education of the Czech Republic, and by the European Union Research Training Network COMBSTRU (Combinatorial Structure of Intractable Problems).}

\subjclass[2000]{05C15 (coloring of graphs and hypergraphs), 11N99 (multiplicative number theory)}

\begin{abstract} 
Let \F\ be a family of connected bipartite graphs, each with at least three vertices. A proper vertex colouring of a graph $G$ with no bichromatic subgraph in \F\ is \emph{\F-free}. The \emph{\F-free chromatic number} $\chi(G,\F)$ of a graph $G$ is the minimum number of colours in an \F-free colouring of $G$. For appropriate choices of \F, several well-known types of colourings fit into this framework, including acyclic colourings, star colourings, and distance-$2$ colourings. This paper studies \F-free colourings of the cartesian product of graphs. 

Let $H$ be the cartesian product of the graphs $G_1,G_2,\dots,G_d$. Our main result establishes an upper bound on the \F-free chromatic number of $H$ in terms of the maximum \F-free chromatic number of the $G_i$ and the following number-theoretic concept. A set $S$ of natural numbers is \emph{$k$-multiplicative Sidon} if $ax=by$ implies $a=b$ and $x=y$ whenever $x,y\in S$ and $1\leq a,b\leq k$. Suppose that $\chi(G_i,\F)\leq k$ and $S$ is a $k$-multiplicative Sidon set of cardinality $d$. We prove that $\chi(H,\F)\leq1+2k\cdot\max S$. We then prove that the maximum density of a $k$-multiplicative Sidon set is $\Theta(1/\log k)$. It follows that $\chi(H,\F)\leq\Oh{dk\log k}$. We illustrate the method with numerous examples, some of which generalise or improve upon existing results in the literature. \end{abstract}

\maketitle

\mySection{Introduction}{Introduction}

\citet{Sabidussi57} proved that the chromatic number of the cartesian product of a set of graphs equals the maximum chromatic number of a graph in the set. No such result is known for more restrictive colourings (such as acyclic, star, and distance-2 colourings). This paper investigates such colourings  of cartesian products under a general model of restriction, in which arbitrary bichromatic subgraphs are excluded. Our study leads to a number-theoretic problem regarding multiplicative Sidon sets that is of independent interest. This problem is then solved using a combination of number-theoretic and graph-theoretic approaches. 

Let $G$ be a graph with vertex set $V(G)$ and edge set $E(G)$. (All graphs considered are undirected, simple, and finite.)\ A \emph{colouring} of $G$ is a function $c:V(G)\rightarrow\Z$ such that $c(v)\ne c(w)$ for every edge $vw\in E(G)$. A colouring $c$ with $|\{c(v):v\in V(G)\}|\leq k$ is a \emph{$k$-colouring}. The \emph{chromatic number} of $G$, denoted by $\chi(G)$, is the minimum integer $k$ for which there is a $k$-colouring of $G$.

Let \F\ be a family of connected bipartite graphs, each with at least three vertices, called a \emph{forbidden family}. A colouring $c$ of a graph $G$ is \emph{\F-free} if it contains no bichromatic subgraph in \F; that is, $|\{c(v):v\in V(H)\}|\geq3$ for every subgraph $H$ of $G$ that is isomorphic to a graph in \F. The \emph{\F-free chromatic number} of $G$, denoted by $\chi(G,\F)$, is the minimum integer $k$ for which there is an \F-free $k$-colouring of $G$. When $\F=\{H\}$ is a singleton, we write \emph{$H$-free} instead of \F-free, and refer to the \emph{$H$-free chromatic number} $\chi(G,H)$. The framework was introduced by \citet{Albertson-EJC04}; an even more general model of restrictive graph colourings is considered by \citet{NesOdM-TreeDepth-EJC06}.

\F-free colourings correspond to many well-studied types of colourings. Let $P_n$ and $C_n$ respectively be the path and cycle on $n$ vertices. Let $\C:=\{C_n:n\text{ even}\}$. Then \C-free colourings are the \emph{acyclic colourings} \citep{KSZ-JGT97, AB-IJM77, B79, Kostochka76, AB-GMJ78, Grunbaum73, BKW-JLMS99, Mitchem-DMJ74, AMR-RSA91, F78, BFKRS-JGT02, Wegner73, BKW-JLMS99, Borodin-DM79, BKRS01, FGR-IPL02, BKRS00, AMS-IJM96, JMV-IPL06,  Berman-CMB78, Wood-DMTCS05, Dvorak-EUJC08, Xu-AC04}. Here each bichromatic subgraph is a forest. By a further restriction we obtain the $P_4$-free colourings, which are called \emph{star colourings}, since each bichromatic subgraph is a collection of disjoint stars \citep{Dvorak-EUJC08, FRR-JGT04, AMS-IJM96, Albertson-EJC04, AMR-RSA91, NesOdM-03, Wood-DMTCS05}. A colouring is $P_3$-free if and only if every pair of vertices at distance at most two receive distinct colours (called a \emph{distance-$2$} colouring). That is, $\chi(G,P_3)=\chi(G^2)$. Here $G^k$ is the $k$-th \emph{power} of $G$, the graph with vertex set $V(G)$, where two vertices are adjacent in $G^k$ whenever they are at distance at most $k$ in $G$. Often motivated by applications in frequency assignment, colourings of graph powers has recently attracted  much attention \citep{FGR-IPL03, AgnHall-SODA04, AH-SJDM03, ADH-DAM03, AGH-CN00, Kral-SJDM04, MS-JCTB05, HM-JGT03}. By definition, 
\begin{equation*}
\chi(G)\,\leq\,\chi(G,\C)\,\leq\,\chi(G,P_4)\,\leq\,\chi(G,P_3)\enspace.
\end{equation*}

Let $G_1$ and $G_2$ be graphs. The \emph{cartesian product} of $G_1$ and $G_2$,
denoted by  \CP{G_1}{G_2}, is the graph with vertex set 
\begin{equation*}
V(\CP{G_1}{G_2})\,:=\,V(G_1)\times V(G_2)\,:=\,\{(a,v):a\in V(G_1),v\in V(G_2)\}
\enspace,
\end{equation*}
where $(a,v)(b,w)$ is an edge of \CP{G_1}{G_2} if and only if  $ab\in E(G_1)$ and $v=w$,  or $a=b$ and $vw\in E(G_2)$.  Assuming isomorphic graphs are equal, the cartesian product is associative, and $\CCP{G_1}{G_2}{G_d}$ is well-defined.  \citet{Sabidussi57} proved that 
\begin{equation*}
\chi(\CCP{G_1}{G_2}{G_d})\,=\,\max\{\chi(G_i):1\leq i\leq d\}\enspace.
\end{equation*}

This paper studies \F-free colourings of cartesian products. The following upper bound on the \F-free chromatic number of a cartesian product is our main result. Here and throughout the paper, $\gamma=0.5772\ldots$ is Euler's constant, and logarithms are base $\e=2.718\ldots$ unless stated otherwise. 

\begin{theorem}
\thmlabel{Chi}
Let \F\ be a forbidden family. Let $G_1,G_2,\dots,G_d$ be graphs, each with \F-free chromatic number $\chi(G_i,\F)\leq k+1$. Then
\begin{equation*}
\chi(\CCP{G_1}{G_2}{G_d},\F)\,\leq\,2k(kd-k+1)+1\enspace.
\end{equation*}
Moreover, for all $\epsilon>0$ and for large $d>d(k,\epsilon)$,
\begin{equation*}
\chi(\CCP{G_1}{G_2}{G_d},\F)\,\leq\,1+\frac{2\,\e^{\gamma}}{1-\epsilon}\,dk\log k\enspace.
\end{equation*}
\end{theorem}

We actually prove a stronger result than \thmref{Chi} that is expressed in terms of `chromatic span'. This concept is introduced in \secref{Span}. The key lemma of the paper, which relates \F-free colourings of a cartesian product to so-called $k$-multiplicative Sidon sets, is proved in \secref{Key}. In \secref{mSets} we study $k$-multiplicative Sidon sets in their own right. The obtained bounds establish our main colouring results. The remaining sections contain numerous examples of the method, some of which generalise or improve upon existing results in the literature. In particular, we consider distance-$2$ colourings in \secref{P3}, acyclic colourings in \secref{Acyclic}, and star colourings in \secref{P4}. Finally, in \secref{Lp1}, we briefly discuss a generalisation of our method for L($p,1$)-labellings.

\mySection{Chromatic Span}{Span}

Let $c$ be a colouring of a graph $G$. The \emph{span} of $c$ is $\max\{|c(v)-c(w)|:vw\in E(G)\}$. (The number of colours is irrelevant.)\ The \emph{chromatic span} of $G$, denoted by $\Lambda(G)$, is the minimum integer $k$ for which there is a colouring of $G$ with span $k$. Note that $\Lambda(G)\leq k$ if and only if there is a homomorphism from $G$ into $P^k_n$ for some $n$.

Let $[a,b]:=[a,a+1,\dots,b]$ and $[b]:=[1,b]$ for all integers $a\leq b$. We can assume that the range of a $k$-colouring is $[k]$. Thus $\Lambda(G)\leq\chi(G)-1$ for every graph $G$. Conversely, given a colouring $c$ of $G$ with span $k$, let $c'(v):=c(v)\bmod{(k+1)}$ for each vertex $v\in V(G)$. Then $c'$ is a ($k+1$)-colouring of $G$. Thus 
\begin{equation*}
\Lambda(G)\,=\,\chi(G)-1\enspace.
\end{equation*}
This might suggest that chromatic span is pointless. Let the \emph{\F-free chromatic span} of a graph $G$, denoted by $\Lambda(G,\F)$, be the minimum integer $k$ for which there is an \F-free colouring of $G$ with span $k$. 

\begin{lemma}
\lemlabel{Span}
Let \F\ be a forbidden family. For every graph $G$, 
\begin{equation*}
\Lambda(G,\F)+1\leq\chi(G,\F)\leq2\cdot\Lambda(G,\F)+1\enspace.
\end{equation*}
\end{lemma}

\begin{proof}
Obviously $\Lambda(G,\F)\leq\chi(G,\F)-1$. To prove that $\chi(G,\F)\leq2\cdot\Lambda(G,\F)+1$, let $c$ be an \F-free colouring of $G$ with span $k:=\Lambda(G,\F)$. For every vertex $v\in V(G)$, let $c'(v):=c(v)\bmod{(2k+1)}$. Clearly $c'$ is a ($2k+1$)-colouring of $G$. For all $i\in[0,2k]$, let $V_i:=\{v\in V(G):c'(v)=i\}$, and for all $j\in\Z$, let $V_{i,j}:=\{v\in V_i:c(v)=j(2k+1)+i\}$. Thus the $V_i$'s are the colour classes of $c'$ and the $V_{i,j}$'s are the colour classes of $c$. For $S,T\subseteq V(G)$ with $S\cap T=\emptyset$, let $G[S,T]$ be the subgraph of $G$ with vertex set $S\cup T$ and edge set $\{vw\in E(G):v\in S,w\in T\}$. Consider two edges $vw,xy\in E(G)$ with $v,x\in V_{i_1,j_1}$, $w\in V_{i_2,j_2}$, and $y\in V_{i_2,j_3}$.  Since $|c(v)-c(w)|\leq k$ and $|c(x)-c(y)|\leq k$, we have $j_2=j_3$. It follows that each bichromatic subgraph of $c'$ is the union of disjoint bichromatic subgraphs of $c$. In particular, $G[V_{i_1},V_{i_2}]=\cup\{G[V_{i_1,j},V_{i_2,j}]:j\in\Z\}$ or $G[V_{i_1},V_{i_2}]=\cup\{G[V_{i_1,j},V_{i_2,j+1}]:j\in\Z\}$. Since each subgraph $G[V_{i_1,j},V_{i_2,j}]$ (or $G[V_{i_1,j},V_{i_2,j+1}]$ in the second case) is \F-free, $c'$ is \F-free and $\chi(G,\F)\leq2\cdot\Lambda(G,\F)+1$.
\end{proof}

\lemref{Span} cannot be improved in general, since it is easily seen that $\Lambda(P_n^k,P_3)=k$ but $\chi(P_n^k,P_3)=2k+1$. Thus chromatic span is of interest when considering \F-free colourings. We prove the following result, which with \lemref{Span}, implies \thmref{Chi}.

\begin{theorem}
\thmlabel{Lambda}
Let \F\ be a forbidden family. Let $G_1,G_2,\dots,G_d$ be graphs, each with \F-free chromatic span $\Lambda(G_i,\F)\leq k$ \paran{which is implied if $\chi(G_i,\F)\leq k+1$}. Then
\begin{align*}
\Lambda(\CCP{G_1}{G_2}{G_d},\F)\,&\leq\,k(kd-k+1)\enspace,
\text{ and}\\
\chi(\CCP{G_1}{G_2}{G_d},\F)\,&\leq\,2k(kd-k+1)+1\enspace.
\end{align*}
Moreover, for all $\epsilon>0$ and for large $d>d(k,\epsilon)$,
\begin{align*}
\Lambda(\CCP{G_1}{G_2}{G_d},\F)\,&\leq\,\frac{\e^{\gamma}}{1-\epsilon}\,dk\log k\text{, and}\\
\chi(\CCP{G_1}{G_2}{G_d},\F)\,&\leq\,1+\frac{2\,\e^{\gamma}}{1-\epsilon}\,dk\log k\enspace.
\end{align*}
\end{theorem}

\mySection{The Key Lemma}{Key}

Our results depend upon the following number-theoretic concept (where $\N:=\{1,2,\dots\}$ and $\NZ:=\N\cup\{0\}$). 

\begin{definition}
Let $k\in\N$. A set $A\subseteq\N$ is \emph{$k$-multiplicative Sidon}\footnote{\citet{Erdos69, Erdos38, Erdos6869} defined a set $A\subseteq\N$ to be \emph{multiplicative Sidon} if $ab=cd$ implies $\{a,b\}=\{c,d\}$ for all $a,b,c,d\in A$; see \citep{Sarkozy01, Ruzsa06, Ruzsa-JNT99}. Additive Sidon sets have been more widely studied; see the classical papers \citep{Sidon32, Singer38, ET41} and the recent survey by \citet{OBryant-EJC04}.} if for all $x,y\in A$ and for all $a,b\in[k]$, we have $ax=by$ implies $a=b$ and $x=y$. For brevity we write \emph{$k$-multiplicative} rather than $k$-multiplicative Sidon.
\end{definition}

Consider a cartesian product $\GG:=\CCP{G_1}{G_2}{G_d}$ to have vertex set
\begin{equation*}
V(\GG)=\{\vv:\vv=(v_1,v_2,\dots,v_d),v_i\in V(G_i),i\in[d]\}\enspace,
\end{equation*}
where $\vv\ww\in E(\GG)$ if and only if $v_iw_i\in E(G_i)$ for some $i$, and
$v_j=w_j$ for all $j\ne i$; we say that the edge $\vv\ww$ is in \emph{dimension} $i$. 

\begin{lemma}
\lemlabel{Key}
Let \F\ be a forbidden family. Let $G_1,G_2,\dots,G_d$ be graphs, each with \F-free chromatic span $\Lambda(G_i,\F)\leq k$ \paran{which is implied if $\chi(G_i,\F)\leq k+1$}. Let $S:=\{s_1,s_2,\dots,s_d\}$ be a $k$-multiplicative set. Then 
\begin{equation*}
\Lambda(\CCP{G_1}{G_2}{G_d},\F)\,\leq\,k\cdot\max S\enspace.
\end{equation*}
\end{lemma}

\begin{proof} 
Let $\GG:=\CCP{G_1}{G_2}{G_d}$. For each $i\in[d]$, let $c_i$ be an \F-free colouring of $G_i$ with span $k$. For each vertex $\vv\in V(\GG)$, let
\begin{equation*}
c(\vv)\,:=\,\sum_{i\in[d]} s_i\cdot c_i(v_i)\enspace.
\end{equation*}
For every edge $\vv\ww\in E(\GG)$ in dimension $i$,
\begin{equation}
\eqnlabel{A}
c(\ww)-c(\vv)\,=\,
\bracket{\sum_{j\in[d]} s_j\cdot c_j(w_j)}-
\bracket{\sum_{j\in[d]} s_j\cdot c_j(v_j)}
\,=\,s_i\big(c_i(w_i)-c_i(v_i)\big)\enspace.
\end{equation}
Since $1\leq|c_i(w_i)-c_i(v_i)|\leq k$ and $s_i\geq1$,  $c$ is a colouring of \GG\ with span $k\cdot\max S$. 

Suppose, for the sake of contradiction, that $c$ is not \F-free. That is, there is a bichromatic subgraph $H$ of \GG\ that is isomorphic to some graph in \F. First suppose that all the edges of $H$ have the same dimension $i$. By \Eqnref{A}, and since $H$ is connected, the edges $\{v_iw_i:\vv\ww\in E(H)\}$ induce a bichromatic subgraph of $G_i$ that is isomorphic to a graph in \F, which is a contradiction. Thus not all the edges of $H$ are in the same dimension. Since $H$ is connected and has at least three vertices, $H$ has two edges $\vv\xx$ and $\ww\xx$ with a common endpoint that are in distinct dimensions. Say $\vv\xx$ is in dimension $i$ and $\ww\xx$ is in dimension $j\ne i$. Since $H$ is bichromatic, $c(\vv)-c(\xx)=c(\ww)-c(\xx)$. By \Eqnref{A},
\begin{equation*}
s_i\big(c_i(v_i)-c_i(x_i)\big)\,=\,s_j\big(c_j(w_j)-c_j(x_j)\big)\enspace.
\end{equation*}
Since $c_i$ has span $k$, we have $1\leq|c_i(v_i)-c_i(x_i)|\leq k$ and $1\leq|c_j(w_j)-c_j(x_j)|\leq k$, which implies that $S$ is not $k$-multiplicative. This contradiction proves that $c$ is an \F-free colouring of \GG. \end{proof}

\mySection{$k$-Multiplicative Sidon Sets}{mSets}

Motivated by \lemref{Key}, in this section we study $k$-multiplicative sets. We measure the `size' of a $k$-multiplicative set by its density. The \emph{density} of $A\subseteq\N$ is 
\begin{equation*}
\Density(A):=\lim_{n \rightarrow \infty} \frac{|A\cap[n]|}{n}
\end{equation*}
if the limit exists (otherwise the density is undefined). We say $A\subseteq\N$ is \emph{$p$-periodic} if $x\in A$ if and only if $x+p\in A$ for all $x\in\N$. Observe that if $A$ is $p$-periodic then 
\begin{equation}
\eqnlabel{Periodic}
\Density(A)=\frac{|A\cap[p]|}{p}.
\end{equation}

The following theorem is our main result regarding $k$-multiplicative sets.

\begin{theorem}
\thmlabel{MultSet}
For all $k\in\N$, the maximum density of a $k$-multiplicative set is 
\begin{equation*}
\Theta\bracket{\frac{1}{\log k}}.
\end{equation*}
\end{theorem}

We start with a naive construction of a $k$-multiplicative set.

\begin{lemma} \lemlabel{BasicConstruction} For all $k\in\N$, the set $R_k:=\{x\in\N:x\equiv1\pmod{k}\}$ is $k$-multiplicative and has density $\Density(R_k)=1/k$. \end{lemma}

\begin{proof} Suppose that $ax=by$ for some $x,y\in R_k$ and $a,b\in[k]$. Then $x=pk+1$ and $y=qk+1$ for some $p,q\in\N$. Thus $(ap-bq)k=b-a$. Since $|b-a|\leq k-1$, we have $a=b$ and $ap=bq$. Thus $p=q$ and $x=y$. That is,  $R_k$ is $k$-multiplicative. Since $R_k$ is $k$-periodic, $\Density(R_k)=|R_k\cap[k]|/k=1/k$ by \Eqnref{Periodic}. \end{proof}

The lower and upper bounds in \thmref{MultSet} are proved in \twothmref{Sk}{MultSetUpperBound}, respectively. 
Fix $k\in\N$. Let $\PP_k:=\{p_1,p_2,\dots,p_\ell\}$ be the set of primes in $[k]$. Let 
\begin{equation*}
\PrimeProduct{k}\,:=\,\prod_{i\in[\ell]}p_i\enspace.
\end{equation*}
Every $x\in\N$ can be uniquely represented as
\begin{equation*}
x\,=\,\beta_*(x)\prod_{i\in[\ell]}p_i^{\beta_i(x)}\enspace,
\end{equation*}
where $\beta_i(x)\in\NZ$ and $\beta_*(x)$ is not divisible by $p_i$ for all $i\in[\ell]$. That is, $\gcd(\beta_*(x),\PrimeProduct{k})=1$. Let $\beta(x)$ be the vector $(\beta_1(x),\beta_2(x),\dots,\beta_\ell(x))$. For all $x,y\in\N$,
\begin{equation}
\eqnlabel{Vectors}
\beta(x\cdot y)\,=\,\beta(x)+\beta(y)
\text{ and }
\beta_*(x\cdot y)\,=\,\beta_*(x)\cdot\beta_*(y)
\enspace.
\end{equation}

\begin{lemma}
\lemlabel{Basic}
For all $k\in\N$, if $ax=by$ for some $a,b\in[k]$ and $x,y\in\N$, then $\beta_*(x)=\beta_*(y)$.
\end{lemma}

\begin{proof}
By \Eqnref{Vectors}, we have $\beta_*(a)\cdot\beta_*(x)=\beta_*(b)\cdot\beta_*(y)$. Since $a,b\leq k$, we have $\beta_*(a)=\beta_*(b)=1$. Thus $\beta_*(x)=\beta_*(y)$.
\end{proof}

\begin{theorem}
\thmlabel{Sk}
For all $k\in\N$, the set $S_k:=\{s\in\N:\gcd(s,\PrimeProduct{k})=1\}$ is $k$-multiplicative and has density
\begin{equation*}
\Density(S_k)
\,=\,\prod_{i\in[\ell]}\bracket{1-\frac{1}{p_i}}
\,\sim\,
\frac{\e^{-\gamma}}{\log k}.
\end{equation*}
\end{theorem}

\begin{proof}
Suppose that $ax=by$ for some $a,b\in[k]$ and $x,y\in S_k$. Thus $\beta_*(x)=\beta_*(y)$ by \lemref{Basic}. Since $\gcd(x,\PrimeProduct{k})=\gcd(y,\PrimeProduct{k})=1$, we have $\beta_i(x)=\beta_i(y)=0$ for all $i\in[\ell]$. Hence $x=y$, which implies that $a=b$, and $S_k$ is $k$-multiplicative. 

It remains to compute the density of $S_k$. Let $\varphi$ be Euler's totient function, $\varphi(x):=|\{y\in[x]:\gcd(x,y)=1\}|$. If $q_1,q_2,\dots,q_r$ are the prime factors of $x$ (with repetition), then 
\begin{equation*}
\varphi(x)\,=\,x\prod_{i\in[r]}\bracket{1-\frac{1}{q_i}}\enspace.
\end{equation*}
Observe that $S_k$ is $\PrimeProduct{k}$-periodic. By \Eqnref{Periodic},
\begin{equation*}
\Density(S_k)
\,=\,\frac{|S_k\cap[\PrimeProduct{k}]|}{\PrimeProduct{k}}
\,=\,\frac{\varphi(\PrimeProduct{k})}{\PrimeProduct{k}}
\,=\,\prod_{i\in[\ell]}\bracket{1-\frac{1}{p_i}}\enspace.
\end{equation*}
By Mertens' Theorem (see \citep{HW79}), $\Density(S_k)\sim \e^{-\gamma}/\log k$; see \tabref{TableS}. 
\end{proof}

The following corollary is a straightforward consequence of \thmref{Sk}.

\begin{corollary}
\corlabel{Sk}
For all $k\in\N$, $\epsilon>0$, and sufficiently large $n>n(k,\epsilon)$,
\begin{equation*}
\frac{(1-\epsilon)n}{\e^{\gamma}\log k}
\,\leq\,
|S_k\cap[n]|
\,\leq\,
\frac{(1+\epsilon)n}{\e^{\gamma}\log k}\enspace.
\end{equation*}
\end{corollary}

\begin{table}[H]
\caption{\tablabel{TableS}
The first $15$ elements of the set $S_k$ for each $k\leq30$.}
\begin{tabular}{clc}
\hline
$k$	& $S_k$	& density\\\hline
$2$
& $\{1,3,5,7,9,11,13,15,17,19,21,23,25,27,29,\dots\}$
& $1/2$\\
$3,4$
& $\{1,5,7,11,13,17,19,23,25,29,31,35,37,41,43,\dots\}$
& $1/3$\\
$5,6$
& $\{1,7,11,13,17,19,23,29,31,37,41,43,47,49,53,\dots\}$
& $4/15$\\
$7,\ldots,10$
& $\{1,11,13,17,19,23,29,31,37,41,43,47,53,59,61,\dots\}$
& $8/35$\\
$11,12$
& $\{1,13,17,19,23,29,31,37,41,43,47,53,59,61,67,\dots\}$
& $16/77$\\
$13,\ldots,16$
& $\{1,17,19,23,29,31,37,41,43,47,53,59,61,67,71,\dots\}$
& $192/1001$\\
$17,18$
& $\{1,19,23,29,31,37,41,43,47,53,59,61,67,71,73,\dots\}$
& $3072/17017$\\
$19,\ldots,22$
& $\{1,23,29,31,37,41,43,47,53,59,61,67,71,73,79,\dots\}$
& $55296/323323$\\
$23,\ldots,28$
& $\{1,29,31,37,41,43,47,53,59,61,67,71,73,79,83,\dots\}$
& $110592/676039$\\
$29,30$
& $\{1,31,37,41,43,47,53,59,61,67,71,73,79,83,89,\dots\}$
& $442368/2800733$\\
\hline
\end{tabular}
\end{table}

We can now prove \thmref{Lambda}.

\begin{proof}[Proof of \thmref{Lambda}] 
\lemref{BasicConstruction} implies that  $R:=\{ik+1:i\in[0,d-1]\}$ is $k$-multiplicative. Since $|R|=d$ and $\max R=dk-k+1$, by using $R$ as a $k$-multiplicative set in \lemref{Key}, we have $\Lambda(\CCP{G_1}{G_2}{G_d},\F)\leq k(dk-k+1)$. This proves the first part of the theorem.  

Let $n$ be the minimum integer such that $|S_k\cap[n]|\geq d$. By \corref{Sk}, for $d>d(k,\epsilon)$, 
\begin{equation*}
\max\{S_k\cap[n]\}\,\leq\,n\,\leq\,\frac{\e^{\gamma}}{1-\epsilon}\,d\log k
\enspace.
\end{equation*}
Using $S_k\cap[n]$ as a $k$-multiplicative set in \lemref{Key}, we have 
\begin{equation*}
\Lambda(\CCP{G_1}{G_2}{G_d},\F)\,\leq\,\frac{\e^{\gamma}}{1-\epsilon}\,dk\log k\enspace.
\end{equation*}
The final claim in \thmref{Lambda} follows from \lemref{Span}.
\end{proof}

\subsection{Proof of Optimality}

We now prove that the lower bound in \thmref{Sk} is asymptotically optimal, which in turn completes the proof of \thmref{MultSet}.

\begin{theorem}
\thmlabel{MultSetUpperBound}
For all $k\in\N$, $\epsilon>0$, and sufficiently large $n>n(k,\epsilon)$, 
every $k$-multiplicative set $A\subseteq[n]$ satisfies
\begin{equation*}
|A|
\,\leq\,\frac{(2+\epsilon)n}{\e^{\gamma}\log k}+\frac{2n}{\sqrt[4]{k}}
\,=\,(2+o(1))|S_k\cap[n]|
\,=\,\frac{(2+o(1))n}{\e^{\gamma}\log k}
\enspace.
\end{equation*}
\end{theorem}

To prove \thmref{MultSetUpperBound}, we model $k$-multiplicative sets using graphs. Let $G_{n,k}$ be the graph with vertex set $V(G_{n,k}):=[n]$, where $xy\in E(G_{n,k})$ whenever $ax=by$ for some $a,b\in[k]$. Observe that a set $A\subseteq[n]$ is $k$-multiplicative if and only if $A$ is an independent set of $G_{n,k}$. For each $s\in S_k\cap[n]$, let $G_{n,k,s}$ be the subgraph of $G_{n,k}$ induced by $X_{n,k,s}:=\{x\in[n]:\beta_*(x)=s\}$.

\begin{lemma}
\lemlabel{ConnectedComponents}
The connected components of $G_{n,k}$ are $\{G_{n,k,s}:s\in S_k\cap[n]\}$.
\end{lemma}

\begin{proof} 
If $xy\in E(G_{n,k})$, then $\beta_*(x)=\beta_*(y)$ by \lemref{Basic}, which implies that $x,y\in X_{n,k,s}$ for some $s\in S_k\cap[n]$. Thus distinct sets $X_{n,k,s}$ and $X_{n,k,t}$ are not joined by an edge of $G_{n,k}$. It remains to prove that each subgraph $G_{n,k,s}$ is connected. For each pair of vertices $x,y\in X_{n,k,s}$, let 
\begin{equation*}
f(x,y):=\sum_{i\in[\ell]}|\beta_i(x)-\beta_i(y)|.
\end{equation*} 
We claim that $x$ and $y$ are connected by a path of $f(x,y)$ edges in $G_{n,k,s}$. The proof is by induction on $f(x,y)$. If $f(x,y)=0$ then $x=y$ (since $\beta_*(x)=\beta_*(y)=s$) and we are done. Say $f(x,y)>0$. Without loss of generality, $\beta_i(x)<\beta_i(y)$ for some $i$. Let $z:=p_ix$. Then $z\in X_{n,k,s}$ and $xz$ is an edge of $G_{n,k,s}$. Moreover, $\beta_i(z)=\beta_i(x)+1$, which implies that $f(z,y)=f(x,y)-1$. By induction, there is a path of $f(z,y)$ edges from $z$ to $y$. Thus there is a path of $f(z,y)+1=f(x,y)$ edges from $x$ to $y$. 
\end{proof}

\begin{lemma}
\lemlabel{StartClique} 
Let $G_{n,k,s}$ be a connected component of $G_{n,k}$ with $r$ vertices. Then the $\min\{k,r\}$ smallest elements of $X_{n,k,s}$ are $\{s,2s,3s,\dots,\min\{k,r\}\cdot s\}$, and they form a clique of $G_{n,k,s}$.
\end{lemma}

\begin{proof}
Every element of $X_{n,k,s}$ is a multiple of $s$ and is at least $s$. Now $is\in X_{n,k,s}$ for each $i\in[\min\{k,r\}]$. Thus the $\min\{k,r\}$ smallest elements of $X_{n,k,s}$ are $\{s,2s,3s,\dots,\min\{k,r\}\cdot s\}$, which clearly form a clique of $G_{n,k,s}$.
\end{proof}

For all $x\in[n]$, let $N_k(x)$ be the closed neighbourhood of $x$ in $G_{n,k}$. That is, $y\in N_k(x)$ if and only if $y\in[n]$ and $ay=bx$ for some $a,b\in[k]$.

\begin{lemma}
\lemlabel{Neighbourhood}
Let $G_{n,k,s}$ be a connected component of $G_{n,k}$ with at least $k$ vertices. 
Then $|N_k(x)|\geq\floor{\sqrt{k}}$ for every $x\in X_{n,k,s}$.
\end{lemma}

\begin{proof} 
By \lemref{StartClique}, the $k$ smallest elements of $X_{n,k,s}$ are $\{s,2s,3s,\dots,ks\}$, and they form a clique of $G_{n,k,s}$. In particular, $ks\leq n$. 

Case (a). $x\leq\sqrt{k}s$: For each $a\in[\floor{\sqrt{k}}]$, we have $ax\leq ks\leq n$. Thus $ax\in N_k(x)$ and $|N_k(x)|\geq\floor{\sqrt{k}}$.

Case (b). $x>\sqrt{k}s$: First suppose that there is a prime $p$ that divides $x$ and $\sqrt{k}\leq p\leq k$. Then $\frac{ax}{p}\in[x]$ for each $a\in[p]$. Thus $\frac{ax}{p}\in N_k(x)$ and $|N_k(x)|\geq p\geq \sqrt{k}$. Now suppose that there is no prime divisor $p$ of $x$ with $\sqrt{k}\leq p\leq k$. Let $p_1\leq p_2\leq \dots\leq p_t$ be the prime factors of $x$ with duplication. Since $x>\sqrt{k}$, for some $\ell\in[t]$, the integer $q:=\prod_{i\in[\ell]}p_i$ divides $x$ and $\sqrt{k}\leq q\leq k$. Thus $\frac{ax}{q}\in[x]$ for each $a\in[q]$. Thus $\frac{ax}{q}\in N_k(x)$ and $|N_k(x)|\geq q\geq\sqrt{k}$. \end{proof}


\begin{proof}[Proof of \thmref{MultSetUpperBound}]
Let $k':=\floor{\sqrt{k}}$ and $k'':=\floor{\sqrt{k'}}$. Note that $k''\geq1$ and $k''>\sqrt[4]{k}/2$. We proceed by studying the size of $A$ within each connected component of the graph $G_{n,k'}$. That is, we consider $A$ as the union of the disjoint sets $\{A\cap X_{n,k',s}:s\in S_{k'}\cap[n]\}$.

First consider $s\in S_{k'}\cap[n]$ for which  $|X_{n,k',s}|\leq k'$. By \lemref{StartClique}, $X_{n,k',s}$ is a clique of $G_{n,k'}$. Since $A$ is $k$-multiplicative, $A$ is $k'$-multiplicative, and $A$ is an independent set of $G_{n,k'}$. Thus $|A\cap X_{n,k',s}|\leq 1$. The set $S_{k'}\cap[n]$ has exactly one element in $X_{n,k',s}$. Thus
$|\cup\{A\cap X_{n,k',s}:s\in S_{k'}\cap[n],|X_{n,k',s}|\leq k'\}|\leq|S_{k'}\cap[n]|$. By \corref{Sk}, 
\begin{equation}
\CARD{\,\bigcup\{A\cap X_{n,k',s}:s\in S_{k'}\cap[n],|X_{n,k',s}|\leq k'\}\,}
\,\leq\,
\frac{(1+\epsilon)n}{\e^{\gamma}\,\log k'}
\,\leq\,
\frac{(2+\epsilon)n}{\e^{\gamma}\log k}.
\eqnlabel{Small}
\end{equation}

Now consider $s\in S_{k'}\cap[n]$ for which $|X_{n,k',s}|>k'$. 
We claim that $N_{k'}(x)\cap N_{k'}(y)=\emptyset$ for distinct $x,y\in A$. 
Suppose that $z\in N_{k'}(x)\cap N_{k'}(y)$ for some $x,y\in A$. Then $a_1x=b_1z$ and $a_2y=b_2z$ for some $a_1,a_2,b_1,b_2\in[k']$. Thus $z=a_1x/b_1=a_2y/b_2$ and $(a_1b_2)x=(a_2b_1)y$. Since $a_1b_2,a_2b_1\in[k]$ and $A$ is $k$-multiplicative, $x=y$. This proves the claim. Now $N_{k'}(x)\subseteq X_{n,k',s}$ for each $x\in X_{n,k',s}$ by \lemref{ConnectedComponents}, and $|N_{k'}(x)|\geq k''$ by \lemref{Neighbourhood}. Thus $|A\cap X_{n,k',s}|\cdot k''\leq|X_{n,k',s}|$, and
\begin{equation}
\CARD{\,\bigcup\{A\cap X_{n,k',s}:s\in S_{k'}\cap[n],|X_{n,k',s}|>k'\}\,}
\,\leq\,
\frac{n}{k''}
\,<\,
\frac{2n}{\sqrt[4]{k}}
\enspace.
\eqnlabel{Big}
\end{equation}

\corref{Sk} and \twoEqnref{Small}{Big} imply that
\begin{align*}
|A|\,
\,\leq\,\frac{(2+\epsilon)n}{\e^{\gamma}\log k}+\frac{2n}{\sqrt[4]{k}}
\,\leq\,\frac{(2+o(1))n}{\e^{\gamma}\log k}
\,=\,(2+o(1))|S_k\cap[n]|\enspace.
\end{align*}
\end{proof}

\subsection{An Improved Construction}

While $S_k\cap[n]$ is a $k$-multiplicative set whose cardinality is within a constant factor of optimal, larger $k$-multiplicative sets in $[n]$ can be constructed. Recall that $\PP_k=\{p_1,p_2,\dots,p_\ell\}$ is the set of primes in $[k]$. Let $\alpha_i:=\floor{\log_{p_i} k}+1$ for each $p_i\in\PP_k$. Define
\begin{align*}
T_k\,&:=\,\{x\in\N\,:\,\beta_i(x)\equiv0\pmod{\alpha_i},\,i\in[\ell]\}\enspace.
\end{align*}

\begin{lemma}
\lemlabel{Tset}
For each $k\in\N$, the set $T_k$ is $k$-multiplicative.
\end{lemma}

\begin{proof}
Suppose that $ax=by$ for some $a,b\in[k]$ and $x,y\in T_k$. By \Eqnref{Vectors},
\begin{equation*}
\beta_i(a)+\beta_i(x)=\beta_i(b)+\beta_i(y)
\end{equation*}
for all $i\in[\ell]$. Now $\beta_i(x)\equiv\beta_i(y)\equiv0\pmod{\alpha_i}$ since $x,y\in T_k$. Thus $\beta_i(a)\equiv\beta_i(b)\pmod{\alpha_i}$. 
Now $p_i^{\beta_i(a)}\leq a\leq k$. Thus $\beta_i(a)\leq\floor{\log_{p_i}k}=\alpha_i-1$. Similarly $\beta_i(b)\leq\alpha_i-1$. Hence $\beta_i(a)=\beta_i(b)$ for all $i\in[\ell]$. Thus $a=b$ and $x=y$. Therefore $T_k$ is $k$-multiplicative.
\end{proof}

We now set out to determine the density of $T_k$. Observe that $S_k=\{x\in\N:\beta_i(x)=0,i\in[\ell]\}\subset T_k$. Thus (if it exists) the density of $T_k$ is at least that of $S_k$. 

Consider $A,B\subseteq\N$ with $A\cap B=\emptyset$. If $\Density(A)$ and $\Density(B)$ exist, then $\Density(A\cup B)=\Density(A)+\Density(B)$. The following lemma extends this idea to an infinite union, where 
\begin{equation*}
\UpperDensity(A)\,:=\,\sup_{n \rightarrow \infty} \frac{|A\cap[n]|}{n}
\enspace.
\end{equation*}

\begin{lemma}
\lemlabel{Technical}
Let $A_1,A_2,\ldots\subseteq\N$ such that $A_i\cap A_j=\emptyset$ whenever $i\ne j$. Suppose that for each $i\in\N$, $\Density(A_i)$ exists and $\UpperDensity(A_i) \leq c\cdot\Density(A_i)$ for some constant $c\geq1$. Let $A:=\bigcup_iA_i$. Then $\Density(A)=\sum_i \Density(A_i)$.
\end{lemma}

\begin{proof}
Let $\Density := \sum_i \Density(A_i)$.
Let $\epsilon>0$ be an arbitrary positive number. 
Let $r_\epsilon$ be the least integer such that
\begin{equation*}
\sum_{i>r_\epsilon} \Density(A_i)\,<\,\frac{\epsilon}{c}\enspace.
\end{equation*}
Let $n_\epsilon$ be the minimum integer such that for all $n>n_\epsilon$
and for all $i\in[r_\epsilon]$,
\begin{equation*}
\CARD{ \frac{|A_i \cap [n]|}{n} - \Density(A_i) }
\,<\,
\frac{\epsilon}{r_\epsilon}\enspace.
\end{equation*}
Let $n>n_\epsilon$, $X := A \cap [n]$, $X_i := X \cap A_i$ and $X^* := \cup\{X_i:i>r_\epsilon\}$. We have  $|X_i| < c\cdot \Density(A_i) n$. Thus
\begin{equation*}
|X^*|\,<\,cn \sum_{i>r_\epsilon} \Density(A_i)\,<\,\epsilon n\enspace.
\end{equation*}
Therefore 
\begin{align*}
	\left | \frac{|X|}{n} - \Density \right | &=
	\left | \BRACKET{\sum_{i\in[r_\epsilon]} \frac{|X_i|}{n} - \Density(A_i) } + 
	\frac{|X^*|}{n} - \sum_{i>r_\epsilon} \Density(A_i) \right | \\
	&< \sum_{i\in[r_\epsilon]} \left | \frac{|X_i|}{n} - \Density(A_i) \right |
	+ \frac{|X^*|}{n} + \sum_{i>r_\epsilon} \Density(A_i) \\
	&< r_\epsilon \frac{\epsilon}{r_\epsilon} + \frac{\epsilon n}{n} + \frac{\epsilon}{c}
	\;<\; \epsilon \bracket{2+\frac{1}{c}} \;<\; 3\epsilon\enspace.
\end{align*}
This proves that $\Density(A) = \Density$.
\end{proof}

\begin{theorem}
\thmlabel{Tk}
The set $T_k$ is $k$-multiplicative with density
\begin{equation*}
\Density(T_k)
\,=\,\Density(S_k) \prod_{i\in[\ell]} \bracket{1+\frac{1}{p_i^{\alpha_i}-1}}
\,=\, \prod_{i\in[\ell]}\bracket{1-\frac{1}{p_i}}\bracket{1+\frac{1}{p_i^{\alpha_i}-1}}\enspace.
\end{equation*}
\end{theorem}

\begin{proof}
For all $A\subseteq\N$ and $t\in\N$, let $t\cdot A:=\{ta:a\in A\}$. 
If $\Density(A)$ exists then 
\begin{equation}
\eqnlabel{MultiplyDivide}
\Density(t\cdot A)\,=\,\frac{\Density(A)}{t}\enspace.
\end{equation}


Now, for all $v\in\NZ^\ell$, let 
\begin{equation*}
S_k^v \,:=\, \bracket{\prod_{i\in[\ell]} p_i^{v_i\alpha_i}}\cdot S_k\enspace.
\end{equation*}
Note that $S_k^v \cap S_k^w = \emptyset$ for distinct $v,w \in \NZ^\ell$. 
For all $v \in \NZ^\ell$ we have 
$\frac{\UpperDensity(S_k^v)}{\Density(S_k^v)}= 
 \frac{\UpperDensity(S_k)}{\Density(S_k)}$.
Now $T_k = \bigcup\{S_k^v:v\in \NZ^\ell\}$. By \lemref{Technical},
\begin{equation*}
\Density(T_k) 
	\,=\, \sum_{v \in \NZ^\ell} \Density(S_k^v)\enspace. 
\end{equation*}
By \Eqnref{MultiplyDivide} with $A=S_k$ and $t=\prod_i p_i^{v_i\alpha_i}$,
\begin{equation*}
\Density(T_k) 
	\,=\, \sum_{v \in \NZ^\ell} \Density(S_k) / \prod_{i\in[\ell]} p_i^{v_i\alpha_i}\enspace.
\end{equation*}
Thus
\begin{equation*}
\Density(T_k) 
	\,=\, \Density(S_k) \sum_{v \in \NZ^\ell} \prod_{i\in[\ell]} p_i^{-v_i\alpha_i} 
	\,=\, \Density(S_k) \prod_{i\in[\ell]} 
		\frac{p_i^{\alpha_i}}{p_i^{\alpha_i}-1}
	\,=\, \Density(S_k) \prod_{i\in[\ell]} 
		\bracket{1+\frac{1}{p_i^{\alpha_i}-1}}\enspace.
\end{equation*}
The result follows by substituting the expression for $\Density(S_k)$ from 
\thmref{Sk}; see \tabref{TableT}.
\end{proof}

\begin{table}[H]
\caption{\tablabel{TableT}
The first $15$ elements of the set $T_k$ for each $k\leq15$.}
\begin{tabular}{clc}
\hline
$k$			& $T_k$	& density\\\hline
$2$ 		& $\{1,3,4,5,7,9,11,12,13,15,16,17,19,20,21,\dots\}$ & $2/3$ \\
$3$ 		& $\{1,4,5,7,9,11,13,16,17,19,20,23,25,28,29,\dots\}$ & $1/2$ \\
$4$ 		& $\{1,5,7,8,9,11,13,17,19,23,25,29,31,35,37,\dots\}$ & $3/7$ \\
$5,6$ 		& $\{1,7,8,9,11,13,17,19,23,25,29,31,37,41,43,\dots\}$ & $5/14$ \\
$7$ 		& $\{1,8,9,11,13,17,19,23,25,29,31,37,41,43,47,\dots\}$ & $5/16$ \\
$8$ 		& $\{1,9,11,13,16,17,19,23,25,29,31,37,41,43,47,\dots\}$ & $7/24$ \\
$9,10$ 		& $\{1,11,13,16,17,19,23,25,27,29,31,37,41,43,47,\dots\}$ & $7/26$ \\
$11,12$		& $\{1,13,16,17,19,23,25,27,29,31,37,41,43,47,49,\dots\}$	& $77/312$\\
$13,14,15$	& $\{1,16,17,19,23,25,27,29,31,37,41,43,47,49,53,\dots\}$ & $11/48$\\
\hline
\end{tabular}
\end{table}


We now show that $\Density(T_k)$ approaches $\Density(S_k)$ for large $k$.

\begin{proposition}
For all $k\in\N$, 
\begin{equation*}
\Density(S_k)\,<\,\Density(T_k)\,=\,c_k\cdot\Density(S_k)\enspace,
\end{equation*}
for some constant $c_k\rightarrow1$ for large $k$.
\end{proposition}

\begin{proof}
By the Prime Number Theorem, $\ell\leq\Oh{k/\log k}$. Thus \\
\EndProof{
\begin{align*}
	c_k 
	&\,=\, \prod_i \bracket{1+\frac{1}{p_i^{\alpha_i}-1}}
	\,<\, \prod_i \bracket{1+\frac{1}{k-1}}
	\,\leq\, \bracket{1+\frac{1}{k-1}}^{\Oh{k/\log k}} \\
	& \,\leq\, \exp( \Oh{1/\log k} )\,\rightarrow\, 1\enspace.
\end{align*}}
\end{proof}


The case $k=2$ was previously studied by \citet{Tamura96} and \citet{AABBJPS-DM95}. Observe that $T_2=\{2^{2i}(2j+1):i,j\in\NZ\}$. \thmref{Tk} with $k=2$ was proved by  \citet{AABBJPS-DM95}, who also proved that $T_2$ has the maximum density out of all $2$-multiplicative sets. Interesting relationships with the Thue-Morse sequence were also discovered.

\begin{proposition}[\citep{AABBJPS-DM95}]
\proplabel{GammaTwo}
The set $T_2$ is $2$-multiplicative and has density $2/3$. For all $d\in\N$, the $d$-th smallest element of $T_2$ is at most $3d/2+\Oh{\log d}$. 
\end{proposition}



\begin{theorem}
\thmlabel{LambdaTwo}
Let \F\ be a forbidden family. Let $G_1,G_2,\dots,G_d$ be graphs, each with $\Lambda(G_i,\F)\leq2$ or $\chi(G_i,\F)\leq3$. 
Let $t$ be the $d$-th smallest element of $T_2$. Then 
\begin{align*}
\Lambda(\CCP{G_1}{G_2}{G_d},\F)\,&\leq\,2t\,\leq\,3d+\Oh{\log d}\enspace,
\text{ and}\\
\chi(\CCP{G_1}{G_2}{G_d},\F)\,&\leq\,4t+1\,\leq\,6d+\Oh{\log d}\enspace.
\end{align*}
\end{theorem}

\begin{proof}
By \lemref{Span}, $\chi(G_i,\F)\leq3$ implies $\Lambda(G_i,\F)\leq2$. 
The result follows by applying \lemref{Key} with the $d$ smallest elements of the $2$-multiplicative set $T_2$ from \propref{GammaTwo}. 
\end{proof}

\mySection{$P_3$-free Colourings}{P3}

Recall that a colouring is $P_3$-free if vertices at distance at most two receive distinct colours. Let $\Delta(G)$ be the maximum degree of the graph $G$. Since a vertex and its neighbours receive distinct colours in a $P_3$-free colouring,
\begin{equation}
\eqnlabel{LowerBound}
\chi(G,P_3)\,\geq\,\Delta(G)+1\enspace.
\end{equation}

Let $Q_d:=\CCP{K_2}{K_2}{K_2}$ be the $d$-dimensional hypercube. 
$P_3$-free colourings of $Q_d$ (and more generally, colourings of powers of $Q_d$) have been extensively studied \citep{Zhou-TCS04, Wan-JCO97, KDP-DAM00, NDG-IPL02, Skupien-DMGT95, Ostergard-JCTA04}. \citet{Wan-JCO97} proved that
\begin{equation*}
d+1\leq\chi(Q_d,P_3)\,\leq\,2^{\ceil{\log_2(d+1)}}\,\leq\,2d\enspace.
\end{equation*}
While our methods are not powerful enough to obtain the above upper bound, for grid graphs we have the following result, which was first proved by \citet{FGR-IPL03}.

\begin{example}[\citep{FGR-IPL03}]
\exlabel{Grid}
Every $d$-dimensional grid graph $G:=\CCP{P_{n_1}}{P_{n_2}}{P_{n_d}}$ satisfies 
$\chi(G,P_3)\leq2d+1$, with equality if every $n_i\geq3$.
\end{example}

\begin{proof} 
The lower bound follows from \Eqnref{LowerBound} since $\Delta(G)=2d$ if 
every $n_i\geq3$. Colour the $i$-th vertex in $P_n$ by $i$. We obtain a $P_3$-free colouring of $P_n$ with span $1$. Thus $\Lambda(P_n,P_3)=1$, and the upper bound follows from \thmref{Lambda} with $k=1$. 
\end{proof}

\exref{Grid} highlights the utility of chromatic span. A weaker bound on 
$\chi(G,P_3)$ is obtained if the $P_3$-free chromatic number, $\chi(P_n,P_3)=3$, is used rather than the the $P_3$-free chromatic span, $\Lambda(P_n,P_3)=1$.


\begin{example}
\exlabel{TwoGrid}
Let $G$ be the $d$-dimensional graph $G:=\CCP{P^2_{n_1}}{P^2_{n_2}}{P^2_{n_d}}$. Let $t$ be the $d$-th smallest element of $T_2$. Then 
\begin{equation*}
\chi(G,P_3)\,\leq\,4t+1\,\leq\,6d+\Oh{\log d}\enspace,
\end{equation*}
and if each $n_i\geq5$ then $\chi(G,P_3)\geq4d+1$.
\end{example}

\begin{proof} 
\Eqnref{LowerBound} implies the lower bound since $\Delta(G)=4d$ if 
each $n_i\geq5$. Obviously $\Lambda(P^2_n,P_3)\leq2$. Thus the upper bound follows from \thmref{LambdaTwo}; see
\tabref{Specific}.
\end{proof}

\begin{table}[H]
\caption{\tablabel{Specific}Upper bound on $\chi(G,P_3)$ for $G:=\CCP{P^2_{n_1}}{P^2_{n_2}}{P^2_{n_d}}$
or $G:=\CCP{C_{n_1}}{C_{n_2}}{C_{n_d}}$.}

\begin{tabular}{cllllllllllllllll}
\hline
$d$				& $1$	& $2$	& $3$	& $4$	& $5$	& $6$	& $7$	& $8$	&  $9$	&  $10$	&  $11$	&  $12$	&  $13$	&  $14$	&  $15$	& $\ldots$\\
$\chi(G,P_3)\leq$	& $5$	& $13$	& $17$	& $21$	& $29$	& $37$	& $45$	& $49$	& $53$	& $61$	& $65$	& $69$	& $77$	& $81$	& $85$	& $\ldots$\\
\hline
\end{tabular}
\end{table}


\begin{example}
\exlabel{GenGrid}
Let $G$ be the graph $\CCP{P^k_{n_1}}{P^k_{n_2}}{P^k_{n_d}}$. 
If there exists $n_i,n_j\geq k$ then $\chi(G,P_3)\geq k^2$, 
and if every $n_i\geq2k+1$ then $\chi(G,P_3)\geq2dk+1$. 
As an upper bound, 
\begin{equation*}
\chi(G,P_3)\,\leq\,2k(kd-k+1)+1\enspace.
\end{equation*}
Moreover, for all $\epsilon>0$ and for large $d>d(k,\epsilon)$,
\begin{equation*}
\chi(G,P_3)\,\leq\,1+\frac{2\,\e^{\gamma}}{1-\epsilon}\,dk\log k\enspace.
\end{equation*}
\end{example}

\begin{proof} If $n_i,n_j\geq k$ then $G^2$ contains a $k^2$-vertex clique, and 
$\chi(G,P_3)=\chi(G^2)\geq k^2$. The second lower bound follows from \Eqnref{LowerBound} since $\Delta(G)=2dk$ if every $n_i\geq2k+1$. Obviously $\Lambda(P^k_n,P_3)\leq k$.  Thus the upper bounds follow from \thmref{Lambda}. \end{proof}


\begin{example}
\exlabel{ToroidalGrid}
The $d$-dimensional toroidal grid $G:=\CCP{C_{n_1}}{C_{n_2}}{C_{n_d}}$ satisfies
\begin{equation*}
2d+1\leq\chi(G,P_3)\,\leq\,4t+1\,\leq\,6d+\Oh{\log d}\enspace,
\end{equation*}
where $t$ is the $d$-th smallest element of $T_2$.
\end{example}

\begin{proof}
The lower bound follows from \Eqnref{LowerBound} since $G$ is $2d$-regular.
Say $C_n=(v_1,v_2,\dots,v_n)$. By considering the vertex ordering 
\begin{equation*}
(v_1,v_n;v_2,v_{n-1};\dots;v_i,v_{n-i+1};\dots;v_{\floor{n/2}},v_{\ceil{n/2}})
\enspace,
\end{equation*}
of $C_n$, we see that $C_n\subset P_n^2$.
Thus the upper bound follows from \exref{TwoGrid}.
\end{proof}

\citet{FRR-JGT04} studied $P_4$-free colourings of toroidal grids, and proved that the minimum number of colours is at most $2d^2+d+1$, and at most $2d+1$ in the case that $2d+1$ divides each $n_i$. Thus \exref{ToroidalGrid} gives a linear upper bound on the $P_3$-free chromatic number of toroidal grids, where even for the weaker notion of $P_4$-free colourings, only a quadratic upper bound was previously known.

\begin{example}
\exlabel{GenToroidalGrid}
Let $G$ be the graph $\CCP{C^k_{n_1}}{C^k_{n_2}}{C^k_{n_d}}$.
If $n_i,n_j\geq k$ for some $i\neq j$, then $\chi(G,P_3)\geq k^2$, 
and if every $n_i\geq2k+1$ then $\chi(G,P_3)\geq2dk+1$.
As an upper bound, 
\begin{equation*}
\chi(G,P_3)\leq 4k(2kd-2k+1)+1\enspace.
\end{equation*}
Moreover, for all $\epsilon>0$ and for large $d>d(k,\epsilon)$,
\begin{equation*}
\chi(G,P_3)\,\leq\,1+\frac{4\,\e^{\gamma}}{1-\epsilon}\,d\cdot 
k\log(2k)\enspace.
\end{equation*}
\end{example}

\begin{proof} The lower bounds are the same as in \exref{GenGrid}. As proved in \exref{ToroidalGrid}, $C_n\subset P_n^2$. Thus $C^k_n\subset P_n^{2k}$, and the upper bound follows from \exref{GenGrid}.
\end{proof}

\mySection{Acyclic Colourings}{Acyclic}

Recall that a colouring with no bichromatic cycle is \emph{acyclic}. The following elementary lower bound is well known, where half the average degree of a graph $G$ is denoted by 
\begin{equation*}
\overline{d}(G)\,:=\,\frac{|E(G)|}{|V(G)|}\enspace.
\end{equation*}

\begin{lemma}
\lemlabel{AcyclicLowerBound}
Every graph $G$ \paran{with at least one edge} has acyclic chromatic number 
$\chi(G,\C)>\overline{d}(G)+1$.
\end{lemma}

\begin{proof}
Say $G$ has an acyclic $k$-colouring. Let $n_i$ be the number of vertices in the $i$-th colour class. Let $m_{i,j}$ be the number of edges between the $i$-th and $j$-th colour classes. Thus $m_{i,j}\leq n_i+n_j-1$. Hence
\begin{equation*}
|E(G)|\;=\!\!\sum_{1\leq i<j\leq k}\!\!\!\!m_{i,j}\;
\leq\!\!\sum_{1\leq i<j\leq k}\!\!\!\!n_i+n_j-1\;\;
=\!\sum_{1\leq i\leq k}\!\!(k-1)n_i\;\;-\binom{k}{2}\;
=\;(k-1)|V(G)|-\binom{k}{2}\enspace.
\end{equation*}
Now $k\geq2$ since $G$ has at least one edge. Thus $k\geq(|E(G)|+1)/|V(G)|+1>\overline{d}(G)+1$.
\end{proof}

It is easily seen that a cartesian product satisfies
\begin{equation}
\eqnlabel{DensProd}
\overline{d}(\CCP{G_1}{G_2}{G_d})\,=\,\sum_{i\in[d]}\overline{d}(G_i)\enspace.
\end{equation}

The following theorem, which was proved for paths by \citet{FGR-IPL03}, gives a special case when a $(k+1)$-colouring can be obtained from a colouring with span $k$, rather than the $(2k+1)$-colouring guaranteed by \lemref{Span}.

\begin{proposition}
\proplabel{AcyclicTrees}
For all trees $T_1,T_2,\dots,T_d$, the acyclic chromatic number 
\begin{equation*}
\chi(\CCP{T_1}{T_2}{T_d},\C)\,\leq\,d+1\enspace,
\end{equation*}
with equality if every $|V(T_i)|\geq d$.
\end{proposition}

\begin{proof} 
Let $\GG:=\CCP{T_1}{T_2}{T_d}$. First we prove the lower bound. By \lemref{AcyclicLowerBound} and \Eqnref{DensProd}, and since $|V(T_i)|\geq d$,
\begin{equation*}
\chi(\GG,\C)
\;>\;\overline{d}(\GG)+1
\;=\;1+\sum_{i\in[d]}\frac{|V(T_i)|-1}{|V(T_i)|}
\;=\;d+1-\sum_{i\in[d]}\frac{1}{|V(T_i)|}
\;\geq\;d\enspace.
\end{equation*}
Hence $\chi(\GG,\C)\geq d+1$. 

Now we prove the upper bound. Root each tree $T_i$ at some vertex $r_i$. For each vertex $v\in V(T_i)$, let $c_i(v)$ be the distance between $r_i$ and $v$ in $T_i$. Then $c_i$ is a colouring of $T_i$ with span one. For each vertex $\vv\in V(\GG)$, let
\begin{equation*}
c(\vv)\;:=\;\sum_{i\in[d]} i\cdot c_i(v_i)\enspace.
\end{equation*}
For each edge $\vv\ww\in E(\GG)$ in dimension $i$,
\begin{equation}
\eqnlabel{AA}
c(\ww)-c(\vv)\,=\,
\bracket{\sum_{j=1}^d j\cdot c_j(w_j)}-
\bracket{\sum_{j=1}^d j\cdot c_j(v_j)}
\,=\,i\big(c_i(w_i)-c_i(v_i)\big)
\,=\,\pm i\enspace.
\end{equation}
Thus $c$ is a colouring of \GG\ with span $d$. Let  $c'(\vv):=c(\vv)\bmod{(d+1)}$. Obviously $c'$ is a $(d+1)$-colouring of \GG. We claim that $c'$ is acyclic. 

Consider each edge of $T_i$ to be oriented away from the root $r_i$. Orient each edge $\vv\ww\in E(\GG)$ in dimension $i$ according to the orientation of $v_iw_i$. That is, orient \vv\ to \ww\ so that $c_i(w_i)-c_i(v_i)=1$. Clearly the orientation of \GG\ is acyclic.

Suppose that on the contrary there is a vertex $\vv\in V(\GG)$ that has two incoming edges $\uu\vv$ and $\ww\vv$ for which $c'(\uu)=c'(\ww)$. Thus $c(\uu)\equiv c(\ww)\pmod{(d+1)}$ and 
\begin{equation*}
c(\uu)-c(\vv)\;\equiv\;c(\ww)-c(\vv)\pmod{(d+1)}\enspace.
\end{equation*}
Let $i$ and $j$ be the dimensions of $\uu\vv$ and $\ww\vv$,  respectively. 
By \Eqnref{AA}, 
\begin{equation*}
i(c_i(u_i)-c_i(v_i))\;\equiv\;j(c_j(w_j)-c_j(v_j))\pmod{(d+1)}\enspace.
\end{equation*}
By the orientation of edges, $c_i(u_i)-c_i(v_i)=1$ and $c_j(w_j)-c_j(v_j)=1$. Thus $i\equiv j\pmod{(d+1)}$, which implies that $i=j$. Hence $\uu=\ww$ since $v_i$ has only one incoming edge in $T_i$ (from its parent). Thus every vertex of \GG\ has at most one incoming edge in each bichromatic subgraph $H$ (with respect to the colouring $c'$). Hence $H$ has an acyclic orientation with at most one incoming edge at each vertex. Therefore $H$ is a forest, and $c'$ is the desired acyclic colouring of \GG.
\end{proof}

\mySection{$P_4$-free Colourings}{P4}

Recall that a colouring with no bichromatic $P_4$ is a \emph{star colouring}.

\begin{example}
\exlabel{StarTrees}
For trees $T_1,T_2,\dots,T_d$, the star chromatic number 
\begin{equation*}
\chi(\CCP{T_1}{T_2}{T_d},P_4)\leq 2d+1\enspace.
\end{equation*}
\end{example}

\begin{proof} Root each tree $T_i$ at some vertex $r_i$. For each vertex $v\in V(T_i)$, let $c_i(v)$ be the distance between $r_i$ and $v$ in $T_i$. (This is the same colouring used in \propref{AcyclicTrees}.)\ Obviously $c_i$ is a $P_4$-free colouring of $T_i$ with span one. The result follows from \thmref{Lambda} with $k=1$. Also note that the same lower bound from \propref{AcyclicTrees} applies for the star chromatic number. \end{proof}

\begin{example}
Let $\mathcal{G}$ be a minor-closed graph family that is not the class of all graphs. Then there is a constant $c=c(\mathcal{G})$ such that for all graphs $G_1,G_2,\dots,G_d\in\mathcal{G}$, 
\begin{equation*}
\chi(\CCP{G_1}{G_2}{G_d},P_4)\,\leq\,cd\enspace.
\end{equation*}
\end{example}

\begin{proof} 
\citet{NesOdM-03} proved that there is a constant $c_1$ (bounded by a small quadratic function of the maximum chromatic number of a graph in $\mathcal{G}$) such that  every graph $G\in\mathcal{G}$ has star-chromatic number $\chi(G,P_4)\leq c_1$. By \thmref{Chi}, there is constant $c_2$ (bounded by a small quadratic function of $c_1$) such that $\chi(\CCP{G_1}{G_2}{G_d},P_4)\leq c_2d$.
\end{proof}

\mySection{L($p,1$)-Labellings}{Lp1}

For $p\in\N$, an L($p,1$)-\emph{labelling} of a graph $G$ is a $P_3$-free colouring of $G$ with the additional property that the colours given to adjacent vertices differ by at least $p$. Such colourings arise in frequency assignment problems; see \citep{GriggsYeh, Yeh-DM06, Goncalves-DM08} for example. 
Let $\Lambda_{p,1}(G)$ be the minimum, taken over all L($p,1$)-labellings $c$ of $G$, of
$$\max_{vw\in E(G)}|c(v)-c(w)|\enspace.$$
Let $\lambda_{p,1}(G)$ be the minimum, taken over all L($p,1$)-labellings $c$ of $G$, of
$$\max_{v\in V(G)}c(v)\;-\min_{v\in V(G)}c(v)\enspace.$$
For example,  $\Lambda_{1,1}(G)=\Lambda(G,P_3)$ and $\lambda_{1,1}(G)=\chi(G,P_3)-1$. 
See \citep{FertinRaspaud-DM07, JNSSS-AC00, ST-DAM06, KuoYan, JKV-DAM05} for results on L($p,1$)-labellings of certain cartesian product graphs. 

It is easily seen that \twolemref{Span}{Key} generalise for L($p,1$)-labellings as follows . 

\begin{lemma}
\lemlabel{Lp1Span}
For every graph $G$, 
\begin{equation*}
\Lambda_{p,1}(G)\leq\lambda_{p,1}(G)\leq2\cdot\Lambda_{p,1}(G)+1\enspace.
\end{equation*}
\end{lemma}

\begin{lemma}
\lemlabel{Lp1Key}
Let $G_1,G_2,\dots,G_d$ be graphs, each with $\Lambda_{p,1}(G_i)\leq k$ \paran{which is implied if $\lambda_{p,1}(G_i)\leq k$}. Let $S:=\{s_1,s_2,\dots,s_d\}$ be a $k$-multiplicative set. Then 
\begin{equation*}
\Lambda_{p,1}(\CCP{G_1}{G_2}{G_d})\,\leq\,k\cdot\max S\enspace.
\end{equation*}
\end{lemma}

\twolemref{Lp1Span}{Lp1Key} imply that \thmref{Lambda} generalises as follows.

\begin{theorem}
\thmlabel{Lp1Lambda}
Let $G_1,G_2,\dots,G_d$ be graphs, each with $\Lambda_{p,1}(G_i)\leq k$ \paran{which is implied if $\lambda_{p,1}(G_i)\leq k$}. Then
\begin{align*}
\Lambda_{p,1}(\CCP{G_1}{G_2}{G_d})\,&\leq\,k(kd-k+1)\enspace,
\text{ and}\\
\lambda_{p,1}(\CCP{G_1}{G_2}{G_d})\,&\leq\,2k(kd-k+1)+1\enspace.
\end{align*}
Moreover, for all $\epsilon>0$ and for large $d>d(k,\epsilon)$,
\begin{align*}
\Lambda_{p,1}(\CCP{G_1}{G_2}{G_d})\,&\leq\,\frac{\e^{\gamma}}{1-\epsilon}\,dk\log k\text{, and}\\
\lambda_{p,1}(\CCP{G_1}{G_2}{G_d})\,&\leq\,1+\frac{2\,\e^{\gamma}}{1-\epsilon}\,dk\log k\enspace.
\end{align*}
\end{theorem}

\bigskip\subsection*{Note: } In related recent work, \citet{JMV-IPL06} independently proved \propref{AcyclicTrees}, and \citet{JamMat-GC08} studied acyclic colourings of cartesian products of cliques (Hamming graphs).


\def\soft#1{\leavevmode\setbox0=\hbox{h}\dimen7=\ht0\advance \dimen7
  by-1ex\relax\if t#1\relax\rlap{\raise.6\dimen7
  \hbox{\kern.3ex\char'47}}#1\relax\else\if T#1\relax
  \rlap{\raise.5\dimen7\hbox{\kern1.3ex\char'47}}#1\relax \else\if
  d#1\relax\rlap{\raise.5\dimen7\hbox{\kern.9ex \char'47}}#1\relax\else\if
  D#1\relax\rlap{\raise.5\dimen7 \hbox{\kern1.4ex\char'47}}#1\relax\else\if
  l#1\relax \rlap{\raise.5\dimen7\hbox{\kern.4ex\char'47}}#1\relax \else\if
  L#1\relax\rlap{\raise.5\dimen7\hbox{\kern.7ex
  \char'47}}#1\relax\else\message{accent \string\soft \space #1 not
  defined!}#1\relax\fi\fi\fi\fi\fi\fi} \def\Dbar{\leavevmode\lower.6ex\hbox to
  0pt{\hskip-.23ex \accent"16\hss}D}

\end{document}